\documentclass[12pt,reqno]{amsart}
                                                                                                                                                                                                                                                                                                                                                                                                                                                                                                                                                                                                                                                                        
\usepackage{lipsum}
\usepackage{graphicx}
\usepackage{array}
\usepackage{amssymb}
\usepackage{hyperref}
\usepackage{amsthm}
\usepackage{amsfonts}
\usepackage{amsmath}
\usepackage{bm}
\usepackage{mathrsfs}
\usepackage[all]{xy}
\usepackage{color}
\usepackage{subfigure}
\usepackage{tikz, tikz-cd}
\usepackage{enumerate}
\usepackage{anysize}
\usepackage{amscd}
\usepackage[letterpaper]{geometry}
\usepackage{geometry}
\usepackage{multirow}
\usepackage{array}
\usepackage{booktabs}
\usepackage{musicography}
\usepackage[utf8]{inputenc}
\usepackage{slashed}

\newtheorem{theorem}{Theorem}[section]

\newtheorem{lemma}[theorem]{Lemma}

\theoremstyle{remark}

\theoremstyle{definition}

\numberwithin{equation}{section}
\numberwithin{figure}{section}
\numberwithin{table}{section}
\numberwithin{equation}{section}
\numberwithin{figure}{section}
\numberwithin{table}{section}

\allowdisplaybreaks[3]

\begin{document}

\title{Remarks on Exact G$_{2}$-Structures on Compact Manifolds}

\author{Aaron Kennon} \thanks{The research of the author is partially suported by Simons Foundation Award \#488629 (Morrison)}
\address{Department of Physics UCSB}
\email{akennon@physics.ucsb.edu}

\begin{abstract}
An important open question in G$_{2}$ geometry concerns whether or not a compact seven-manifold can support an exact G$_{2}$-Structure. Given the significance of this question we initiate a study of exact G$_{2}$-Structures on compact manifolds. We focus on exact G$_2$-Structures subject to no additional constraints but we also consider the relationship between the exact condition and other conditions for closed G$_{2}$-Structures such as the Extremally Ricci-Pinched and Laplacian Soliton conditions. 
\end{abstract}

\date{}
\maketitle
\setcounter{tocdepth}{1}
\tableofcontents

\section{Introduction}

Manifolds with G$_{2}$-holonomy metrics are quite distinguished geometric objects. From the perspective of mathematics they are significant because the Lie group G$_{2}$ is the only possible special holonomy group for a simply-connected, irreducible, non-symmetric Riemannian seven-manifold \cite{Berger}. Such manifolds are also Ricci-flat and admit a parallel spinor \cite{Wang} which makes them very interesting as backgrounds for M-Theory compactifications to four-dimensional Minkowski space. \cite{Acharya, AG,AW} \\

In order for a manifold to have G$_{2}$ holonomy it needs to admit a torsion-free G$_{2}$-Structure. Such a structure is equivalent to the existence of a distinguished harmonic three-form $\varphi$ pointwise identified with a model form on $\mathbb{R}^7$. Provided that a compact manifold admits a torsion-free G$_{2}$-Structure the manifold has holonomy precisely equal to G$_{2}$ if and only if the fundamental group is finite. \\

All known constructions of G$_{2}$-holonomy metrics on compact manifolds start from closed G$_{2}$-Structures which are weaker structures than torsion-free G$_{2}$-Structures where the three-form $\varphi$ is only required to be closed and not necessarily co-closed. Not much is known about these closed G$_{2}$-Structures and understanding them better is an important step toward determining when a  manifold admits torsion-free G$_{2}$-Structures. \\

In particular, it is not known if a compact manifold can admit an exact G$_{2}$-Structure. Of course, on any manifold with vanishing third Betti number any closed G$_{2}$-Structure will have to be exact. We can refine this question by asking if compact manifolds subject to additional constraints can admit exact G$_2$-Structures. If exact G$_2$-Structures cannot exist on manifolds with vanishing first Betti number, then that would mean that a manifold which admits a G$_2$-Holonomy metric could not also admit an exact G$_2$-Structure. Similarly, if an exact G$_2$-Structure could not exist on a manifold with nonvanishing third Betti number, then a manifold which admits a torsion-free G$_2$-Structure of any holonomy could not admit an exact G$_2$-Structure. \\

The case of closed G$_{2}$-Structures should be contrasted with that of co-closed G$_{2}$-Structures. In that case we have a parametric h-principle which says that any compact, orientable, spinnable seven-manifold admits a co-closed G$_{2}$-Structure C$^{0}$ close to any given G$_{2}$-Structure \cite{CN}. As a consequence, co-closed G$_2$-Structures are very plentiful even on compact manifolds. It is known that pointwise the closed condition is stronger than the co-closed condition, but it is unclear precisely how restrictive the closed conditon is on compact manifolds. If it turns out that exact G$_{2}$-Structures on compact manifolds do not exist, then compact manifolds with vanishing third Betti number cannot admit closed G$_{2}$-Structures at all, which would imply it is more difficult to satisfy the closed condition on compact manifolds than the co-closed condition. Such a result would be the first nonexistence result for closed G$_{2}$-Structures subject to no additional constraints on compact manifolds that otherwise admit G$_{2}$-Structures.  \\

We begin by introducing background on G$_{2}$-Structures in \S2. We then prove pointwise identities for exact G$_{2}$-Structures and use them to prove results on the two-form underlying an exact G$_{2}$-Structure on a compact manifold in \S3. We finish by considering additional constraints on exact G$_{2}$-Structures on compact manifolds in \S4 and results on the Riemann curvature in \S5.

\section{Background on G$_{2}$-Structures}

Here we include a general introduction to G$_{2}$-Structures on seven-manifolds. The material here may be found in the review by Bryant \cite{Bryant1}, the review by Karigiannis \cite{Karigiannis1}, or the monograph by Joyce \cite{JoyceBook}. \\

A \textit{G$_{2}$-Structure} is a principle subbundle of the frame bundle with structure group G$_{2}$. A seven-manifold admits G$_{2}$-Structures if and only if it is orientable and spin. We may also characterize the existence of a G$_{2}$-Structure in terms of a three-form pointwise identified with the three-form $\varphi_{0}$ on $\mathbb{R}^{7}$.

\begin{equation}
\varphi_0 = dx^{123} - dx^{167} - dx^{527} - dx^{563} + dx^{415} +
dx^{426} + dx^{437}
\end{equation}

Such a three-form is called \textit{positive} and we also refer to such a form as the \textit{G$_2$ three-form}. A positive three-form defines a Riemannian metric and also determines a four-form by taking the Hodge star of the three-form with respect to this metric. Since both the bundle of positive three-forms and that of three-forms on a seven-manifold have rank 35 positivity is an open condition. In other words, if we have a G$_{2}$-Structure and perturb it by a three-form with sufficiently small C$_{0}$ norm, then the resulting form is also a G$_{2}$-Structure. \\

On a manifold with a positive three-form, we may decompose the spaces of differential k-forms into irreducible G$_{2}$ representations. The equations defining the various representations are then all given in terms of the G$_{2}$ three-form $\varphi$ and G$_2$ four-form $\psi$ \cite{Karigiannis3}.

\begin{align}
 \hspace{-100mm} & \Lambda^{1}(T^{*}M) = \Lambda^{1}_{8}(T^{*}M) \\
& \Lambda^{2}(T^{*}M) = \Lambda^{2}_{7}(T^{*}M) \oplus \Lambda^{2}_{14}(T^{*}M) \\
& \Lambda^{3}(T^{*}M) = \Lambda^{3}_{1}(T^{*}M) \oplus \Lambda^{3}_{7}(T^{*}M) \oplus \Lambda^{3}_{27}(T^{*}M) 
\end{align}

\begin{align}
& \Lambda^{2}_{7}(T^{*}M) = \{ \beta \in \Lambda^{2}(T^{*}M) : *(\psi \wedge \beta) = 2\beta \} \\
& \Lambda^{2}_{14}(T^{*}M) =  \{ \beta \in \Lambda^{2}(T^{*}M) : *(\psi \wedge \beta) = -\beta \} \\ 
& \Lambda^{3}_{1}(T^{*}M)= \{  f\varphi \hspace{2mm} \textnormal{for} \hspace{1mm} f \in C^{\infty}(M) \} \\
& \Lambda^{3}_{7}(T^{*}M) =  \{ X \lrcorner \psi \hspace{2mm} \textnormal{for}\hspace{1mm} X \in \Gamma(M) \} \\
& \Lambda^{3}_{27}(T^{*}M) =  \{ \eta \in \Lambda^{3}(T^{*}M): \eta \wedge \varphi = \eta \wedge \psi=0 \} 
\end{align}

We may also write out the representation equations in local coordinates. For instance, for two-forms the equations defining the representations may be expressed easily in terms of contractions between the relevant tensors. 

\begin{align}
& \beta \in  \Lambda^{2}_{7}(T^{*}M) \hspace{2mm} \textnormal{iff} \hspace{2mm} \beta_{ij}\psi_{abij}=4\beta_{ab} \\
& \beta \in  \Lambda^{2}_{14}(T^{*}M) \hspace{2mm} \textnormal{iff} \hspace{2mm} \beta_{ij}\varphi_{ijk}=0
\end{align}

There is an important isomorphism of G$_{2}$ representations $\i_{\varphi}$ which maps the space of symmetric, traceless two-tensors S$_{0}^{2}(M)$ onto $\Lambda_{27}^{3}(T^{\ast}M)$. The mapping may be defined explicitly on decomposable, symmetric, traceless two tensors $\alpha \circ \beta$ \cite{Bryant1}.

\begin{equation}
\i_{\varphi}(\alpha\circ\beta) = 
\alpha\wedge\ast(\beta\wedge\ast\varphi)
+\beta\wedge\ast(\alpha\wedge\ast\varphi)
\end{equation}

As any symmetric traceless two-tensor is built out of decomposable elements we can define this map for an arbitrary such tensor by extending by linearity. We also may explicitly invert the isomorphism $i_{\varphi}$ to get an isomorphism $\j_{\varphi}$ which maps an element $\gamma \in \Lambda_{27}^{3}(T^{\ast}M)$ to a symmeteric traceless two-tensor with vector arguments $v$ and $w$.  
 
\begin{equation}
\j_{\varphi}(\gamma)(v,w) = \ast\bigl((v\lrcorner\varphi)\wedge(w\lrcorner\varphi)\wedge\gamma\bigr)
\end{equation}

In order to ensure a holonomy reduction, we need to find a G$_{2}$-Structure that is torsion-free. The G$_{2}$-Structure being torsion-free is equivalent to the positive three-form being closed and co-closed with respect to the Riemannian metric determined by the three-form. On a compact manifold, a differential form being both closed and co-closed is equivalent to it being harmonic.

\begin{align}
& d\varphi=0 \\
&d\ast\varphi=0
\end{align}

Provided we have a torsion-free G$_{2}$-Structure, the holonomy of the corresponding Levi-Civita connection is known to be contained in G$_{2}$ but it may be a proper subgroup. The possibilities are the trivial group, SU(2), and SU(3). If the holonomy is trivial, then the corresponding connection is flat as a corollary of the Ambrose-Singer theorem \cite{JoyceBook}. Given a torsion-free G$_{2}$-Structure, we have straightforward means to ensure the holonomy is precisely G$_{2}$. On a compact manifold we can do so by ensuring the fundamental group is finite \cite{JoyceBook} and on a simply-connected complete manifold we may ensure that there are no parallel one-forms \cite{BS}. \\

We may generalize the torsion-free condition by studying positive three-forms that are not necessarily harmonic. For an arbitrary G$_{2}$-Structure we may decompose the exterior derivatives of the G$_{2}$ three-form and four-form into irreducible G$_{2}$ representations.

\begin{align}
& d\varphi = \tau_{0}\psi + 3 \tau_{1} \wedge \phi + \ast \tau_{3} \\
& d\psi = 4 \tau_{1} \wedge \psi + \ast \tau_{2}
\end{align}

This decomposition identifies special functions and differential forms $\tau_{0} \in C^{\infty}(M)$, $\tau_{1} \in \Gamma(M)$, $\tau_{2} \in \Lambda^{2}_{14}(T^{*}M)$, and $\tau_{3}\in \Lambda^{3}_{27}(T^{*}M)$. There are then sixteen possibilities for types of G$_{2}$-Structures depending on which of these torsion forms vanish. We can amalgamate all of these torsion forms into a tensor called the torsion tensor. We may express this tensor in terms of the covariant derivative of the three-form and the four-form in coordinates.

\begin{equation}
T_{lm} = 1/24 (\nabla_{l}\varphi_{abc})\psi_{mabc}
\end{equation}

We may also express the torsion tensor in terms of the torsion forms $\tau_{i}$.

\begin{equation}
T_{lm}=\frac{\tau_{0}}{4}g_{lm}-(\tau_{3})_{lm} -(\tau_{1})_{lm} -\frac{1}{2}(\tau_{2})_{lm}
\end{equation}

Here $\tau_{1}$, which properly speaking is a one-form, is identified with a two-form in the $\Lambda^{2}_{7}(T^{*}M)$ representation of G$_{2}$ on two-forms. The curvature of the metric associated to the given G$_{2}$-Structure is completely determined by the torsion tensor. An expression for the Riemann curvature in terms of the torsion tensor is given in the work of Karigiannis \cite{Karigiannis2}, and expressions for the Ricci tensor and scalar curvature of arbitrary G$_2$-Structures in terms of the torsion forms are given by Bryant \cite{Bryant1}. \\

%The Levi-Civita connection is a torsion-free connection but on a manifold with a G$_2$-Structure there is another connection called the canonical G$_2$-connection which has torsion given by the torsion tensor of the G$_2$-Structure. As a torsion-free G$_{2}$-Structure (as its name suggests) has no torsion, in that case this connection corresponds to the Levi-Civita connection. \\

One particularly interesting class of G$_{2}$-Structures corresponds to the case when the positive three-form is closed but not necessarily co-closed. Such G$_{2}$-Structures are called \textit{closed G$_{2}$-Structures} and are sometimes also called \textit{calibrated G$_{2}$-Structures}. Looking at the expression for the exterior derivative of $\varphi$ for a closed G$_{2}$-Structure, it is apparent that every torsion component vanishes except for $\tau_{2}$. As we will be discussing closed G$_{2}$-Structures exclusively, we will abbreviate $\tau_{2}$ by $\tau$. \\

Closed G$_{2}$-Structures are significant as they are used in constructions for manifolds with holonomy G$_{2}$. The idea is to construct a closed G$_{2}$-Structure whose torsion is sufficiently small and then use perturbative analysis to deform the solution to a torison-free G$_{2}$-Structure. The precise criterion needed to perturb these structures are encapsulated in a theorem of Joyce \cite{Joyce1} and a lot of effort in G$_{2}$-Geometry is devoted to the question of when a closed G$_{2}$-Structure can be deformed into a torsion-free G$_2$-Structure. For instance, providing a refinement of our understanding of this issue is a primary motivation of studying the Laplacian flow of closed G$_{2}$-Structures \cite{LW}. \\

The expressions for the torsion tensor and the various cuvature tensors simplify considerably when the G$_{2}$-Structure is closed. The torsion tensor is then directly proportional to the torsion two-form $\tau$ which implies that the torsion tensor is totally antisymmetric for a closed-G$_{2}$-Structure. The Ricci curvature and scalar curvature expressions also simplify considerably. 

\begin{align}
\textnormal{Ric}(g_\varphi) &= 
  \biggl(\frac{1}{4}\,|\tau|^2\biggl)\,g_\varphi - \j_{\varphi}\Bigl( \frac{1}{4}\, d\tau - \frac{1}{8}\,\ast(\tau\wedge\tau \Bigr) \\
\textnormal{Scal}(g_\varphi) &= -\frac{1}{2}|\tau|^2
\end{align}

In particular, the scalar curvature is non-positive and vanishes identically if and only if the G$_2$-Structure is torsion-free. In this case the metric is not only scalar-flat but also Ricci-flat. \\

%Include expression for d$\tau$

Besides these results concerning the torsion and curvature relatively little is known about these structures. There are now several methods of constructing closed G$_{2}$-Structures in the literature. The original approach utilized symmetries to simplify the PDE system $d\varphi=0 \hspace{2mm} $\cite{F1, F2, CF}. The examples have the form $M=\Gamma/G$ where $G$ is a seven-dimensional simply-connected nilpotent or more general solvable Lie group and $\Gamma$ is a co-compact discrete subgroup of $G$. The closed G$_2$-Structures on these spaces come from left-invariant ones on $G$. An advantage of these examples is that the G$_2$-Structure is determined by the three-form on the Lie algebra of G and the closed condition reduces to an algebraic equation. \\

%\textcolor{red}{Discuss compact quotients. Discuss the manifolds more specifically} 

A newer method to constuct closed G$_2$-Structures on compact manifolds is based on the known constructions of G$_2$ holonomy manifolds. A preliminary step for constructing these G$_2$ holonomy metrics is to get a closed G$_2$-Structure on the manifold with small torsion so these methods can be thought of as constructions for closed G$_2$-Structures with more stringent requirements. Of course we may disregard these requirements if the goal is just to get closed G$_2$-Structures. \\

One G$_2$ holonomy construction which has been considered in the context of closed G$_2$-Structures is Joyce's orbifold construction \cite{Joyce1, Joyce2}. In the original version Joyce starts with the standard torsion-free G$_2$-Structure on $T^7$ (with trivial holonomy) and acts on this space with a finite group which preserves the G$_2$-Structure. Then he glues in spaces which resolve these orbifold singularities producing a smooth manifold. On this manifold he uses the original G$_2$-Structures to build a new G$_2$-Structure but this G$_2$ three-form only closed as opposed to harmonic. These structures have arbitrarily small torsion depending on parameters involved in the gluing. For the closed G$_2$-Structure construction the starting manifold is a compact nilmanifold $M=\Gamma/N$ for a particular nilpotent Lie group N \cite{FFKM}. Then for a particular finite group the resulting orbifold may be desingularized and shown to carry a closed G$_2$-Structure. The resulting manifold is interesting because it has nonvanishing first Betti number so it cannot admit a G$_2$ holonomy metric. \\

To date there are no known examples of compact manifolds which admit G$_{2}$-Structures but do not admit closed G$_{2}$-Structures. However, there are a number of geometric constraints that closed G$_{2}$-Structures are known not to be able to satisfy. Several of these have to do with the curvature of the metrics determined by these structures. For instance, the metric corresponding to a closed G$_{2}$-Structure on a compact manifold cannot be Einstein and also cannot be too Ricci-pinched \cite{Bryant1, CI1}. The automorphism group of these structures is also constrained by the topology of the underlying manifold \cite{PR}. We will see that many of these results which are known to hold for closed G$_{2}$-Structures generally may be strengthened when the G$_{2}$-Structure is exact. \\

There are also several classes of exact G$_{2}$-Structures on non-compact manifolds that have recently been produced. Examples on unimodular Lie algebras with vanishing third Betti number are considered by Fernandez, Fino, and Raffero \cite{FFR}, and examples which are closed G$_2$-Eigenforms are considered by Friebert and Salamon \cite{FS}. It is worthwhile to note that there are no compact quotients of the unimodular Lie algebra examples, and it is known that a compact manifold cannot admit closed G$_{2}$-eigenforms \cite{LW}. 

\section{Constraints on the Two-Form Underlying an Exact G$_2$-Structure}

We start by working out a number of projection identities for differential forms onto irreducible G$_2$ components that hold generally when the G$_2$-Structure is closed. For this work as we are mostly interested in exact G$_2$-Structures, which are defined by a two-form, the two-form projection identities are particularly significant. 

\begin{lemma} [Projection identities for derivatives of two-forms] \label{lemma:G$_{2}$ Projections} If (M, $\varphi$) is a manifold with closed G$_{2}$-Structure $\varphi$ , X is a vector field, and $\alpha \in \Lambda^{2}_{14}(T^{*}M)$ then the type decompositions of $d(X\lrcorner\varphi)$ and $d\alpha$ may be given in terms of X, $\alpha$, and $\gamma_{1},\gamma_{2}\in\Lambda^{3}_{27}(T^{*}M)$.
\begin{align}
& d(X\lrcorner\varphi) = -3/7 (d^{\ast}X) \varphi - 1/2(d^{\ast}(X \lrcorner \varphi)) \lrcorner \psi + \gamma_{1}  \\
& d\alpha_{14} = 1/7 \langle\alpha_{14},\tau\rangle\varphi + (1/4d^{\ast}\alpha_{14}) \lrcorner \psi + \gamma_{2} 
%& d^{\ast}(X\lrcorner\varphi) = \\
%&d^{\ast}\beta_{14} =
\end{align}
\end{lemma}

\textit{Proof.} $d(X\lrcorner\varphi)$ and $d\beta_{14}$ are both three-forms so they will take the form $f\varphi + Y \lrcorner \psi +  \gamma$ for some function $f$, vector field $Y$ and three-form $\gamma \in \Lambda^{3}_{27}(T^{*}M)$. In both cases we can isolate the one dimensional component by wedging the relevant three-form with $\psi$. For $d(X\lrcorner\varphi)$ doing so gives a factor of $d^{\ast}X$.

\begin{equation}
d(X\lrcorner\varphi) \wedge \psi = 7f \textnormal{vol}
\end{equation}

We can then use the fact that the exterior derivative is a derivation to rearrange this expression.

\begin{equation}
d(X\lrcorner\varphi \wedge \psi) - X\lrcorner\varphi \wedge \tau \wedge \varphi = 7f \textnormal{vol}
\end{equation}

However, $X\lrcorner\varphi$ and $\tau$ are pointwise orthogonal so the second term vanishes identically. The argument in the differential also simplifies considerably \cite{Karigiannis3}. 

\begin{equation}
X\lrcorner\varphi \wedge \psi = 3\ast X
\end{equation}

We may then take the Hodge star to get the desired expression for $f$. For $d\alpha_{14}$ wedging with $\psi$ gives the pointwise inner product of $\alpha_{14}$ with the torsion two-form $\tau$.

\begin{equation}
d\alpha_{14} \wedge \psi = 7f\textnormal{vol}
\end{equation}

We may then rearrange this expression and use the fact that $\alpha_{14}\wedge\psi=0$ to eliminate the exact term.

\begin{equation}
-\alpha_{14}\wedge\tau\wedge\varphi = 7f\textnormal{vol}
\end{equation}

But the left-hand side is just an expression for the pointwise interior product of $\alpha_{14}$ and $\tau$. The result follows. \\

To isolate the seven-dimensional components of these forms we need to wedge the relevant forms by $\varphi$. For $d(X\lrcorner\varphi)$ we start by using G$_{2}$ identities to get an expression for the Hodge dual of Y \cite{Karigiannis1}.

\begin{equation}
d(X\lrcorner\varphi) \wedge \varphi = -4 \ast Y
\end{equation}

Then we can use the Leibnitz rule to get the answer solely in terms of the vector field X.

\begin{equation}
d(X\lrcorner\varphi \wedge \varphi) = -1/2 d\ast (X \lrcorner \varphi) = -4\ast Y
\end{equation}

We can then take the Hodge star to isolate the vector field Y.

\begin{equation}
-1/2 d^{\ast}(X\lrcorner\varphi) = Y
\end{equation}

For $\alpha_{14}$ we can follow a similar procedure.

\begin{equation}
d\alpha_{14} \wedge \varphi = -4\ast Y
\end{equation}

We can use the equation that defines the 14-dimensional representation to simplify the last equation and get the answer in terms of the co-exterior derivative of $\alpha_{14}$.

\begin{equation}
1/4 (d^{\ast} \alpha_{14}) = Y
\end{equation}

All that is left is the twenty-seven dimensional components of these forms which may be thought of as what is left after subtracting off the components corresponding to the one and seven-dimensional representations. \qed

%\begin{theorem} Projection identities for closed G$_{2}$-Structure three-forms
%\begin{align}
%& d(f\varphi) =\\
%&d^{\ast}(f\varphi) =\\
%& d(X\lrcorner\psi) =  \\
%& d^{\ast}(X\lrcorner\psi) = 
%\end{align}
%\end{theorem}

%\textbf{Four-form projection identities} \\

%\textbf{Lie Derivative Identities} 

Now we apply some of these general projection identities to infer information about the underlying two-form and related differential forms associated to an exact G$_{2}$-Structure. Many of these results require compactness, so we are careful to emphasize when we are utilizing this assumption. \\

Let $\varphi=d\beta$ be an exact G$_2$-Structure. An important point to note is that for such a structure we may add a closed two-form $\alpha$ to $\beta$ and the resulting G$_{2}$-Structure $d\beta$ will be unchanged. In particular, everything determined by the G$_{2}$-Structure such as the torsion or decomposition of forms into irreducible representations of G$_{2}$ will be unchanged under such a transformation. Keeping this gauge freedom in mind, we can fix a particular two-form $\beta$ corresponding to exact G$_{2}$-Structure $d\beta$  and we can decompose $\beta$ into irreducible representations of G$_{2}$.

\begin{equation}
\beta = \beta_{7} + \beta_{14}
\end{equation}

Then if we differentiate we get an expression for the G$_2$-Structure in terms of derivatives of the two-form components. 

\begin{equation}
d\beta = d\beta_{7} + d\beta_{14}
\end{equation}

We can then decompose $d\beta_{7}$ and $d\beta_{14}$ into irreducible G$_{2}$ representations and we get some pointwise constraints on the sums of these components since d$\beta$ is the G$_{2}$-Structure.

\begin{align}
& \pi_{1}(d\beta_{7}) + \pi_{1}(d\beta_{14}) = 1 \\
& \pi_{7}(d\beta_{7}) + \pi_{7}(d\beta_{14}) = 0 \\
& \pi_{27}(d\beta_{7}) + \pi_{27}(d\beta_{14}) = 0 
\end{align}

We can now substitute the equations we found for each of these components into these constraints. Doing so for the one-dimensional component, we get a relation between the pointwise inner product of $\beta_{14}$ and $\tau$ and the co-exterior derivative of $X$.

\begin{equation}
\langle\beta_{14},\tau\rangle -3 d^{\ast}X = 7
\end{equation}

If we integrate this pointwise statement over a compact manifold, the divergence term vanishes and we get an equation for the total volume.

\begin{equation}
\langle\beta_{14},\tau\rangle_{L^2} = 7\textnormal{Vol}
\end{equation}

We could arrive at the same conclusion by considering the $L^2$ norm of the G$_{2}$-Structure $d\beta$. As the G$_{2}$-Structure has pointwise norm equal to seven, the $L^2$ norm gives an expression for the total volume.

\begin{equation}
7\textnormal{Vol} = |d\beta|^2_{L^2}
\end{equation}

Assuming compactness, we can then integrate the $L^2$ norm of $d\beta$ by parts and utilize the fact that $\tau$ is in the fourteen dimensional representation of G$_{2}$ on two-forms to arrive at the result.

\begin{equation}
|d\beta|^2_{L^2} = \langle\beta,\tau\rangle_{L^2} = \langle\beta_{14},\tau\rangle_{L^2}
\end{equation}

As a consequence of this result it is clear that $\beta_{14}$ cannot vanish identically. Later in this section we will prove that $\beta_{7}$ also cannot vanish and we will be able to relate the $L^2$ norms of these two-forms. \\

We now turn to the constraint from $\pi_{7}(d\beta)=0$. From this identity we see that  we get a canonical co-closed two-form built out of the components of the two-form $\beta$ underlying an exact G2-Structure. 

\begin{equation}
d^{\ast} (2\beta_{7}-\beta_{14})=0
\end{equation}

Note that the co-closed form $2\beta_{7} - \beta_{14}$ can't vanish unless its components vanish individually, which would imply that $\beta$ itself vanishes. \\

%\textbf{Consider 27 component if appropriate}

Now we assume compactness and we want to see if $\beta$ even needs to have both components. It turns out that it does, and moreover, we are able to get a constraint relating the $ L^{2}$ norms of these components.

\begin{theorem} [L$^{2}$ norms of $\beta$ components] \label{thm:L2NormsBeta} The L$^{2}$ norms of the irreducible G$_{2}$-components of the two-form $\beta$ underlying an exact G$_2$-Structure on a compact manifold are related up to a factor of two. 
\begin{equation}
2 |\beta_{7}|^{2}_{L^{2}} = |\beta_{14}|^{2}_{L^{2}}
\end{equation}
\end{theorem}

\textit{Proof}. We consider the exterior derivative of $\beta$ wedged against itself three times. We may expand this expression using the product rule. 

\begin{equation}
1/3 \hspace{2mm} d(\beta \wedge \beta \wedge \beta) = \beta \wedge \beta \wedge d\beta
\end{equation}

If we assume compactness then we may integrate the left-hand side of this equation and get zero by citing Stokes' theorem. However, on the right-hand side d$\beta$ is the G$_{2}$-Structure, and we get an identity for a two-form wedged twice against a G$_{2}$-Structure using the representation equations for two-forms and orthogonality of the representations. 

\begin{align}
& \beta \wedge \varphi = 2 \ast \beta_{7} -\ast \beta_{14} \\
& \beta \wedge \beta \wedge \varphi = 2 |\beta_{7}|^{2} - |\beta_{14}|^{2}
\end{align}

But we know that if we integrate this, then that is equal to zero so the result follows. \qed \\

The main implication of this result is that neither irreducible G$_{2}$ component of the two-form $\beta$ can vanish identically if $\beta$ defines an exact G$_{2}$-Structure. Moreover, since $\beta_{7}$ is determined by a one-form and doesn't vanish identically we also have a canonical vector field X associated to an exact G$_{2}$-Structure which cannot vanish. We can then take the derivative of the corresponding one-form to get a two-form dX which we can relate to $\beta$ and $\tau$. \\

%\textcolor{red}{Consider putting constraints here relating dX, $\tau$ and $\beta$. Decomposition of dX into irreducible G$_2$ reps}

%\textcolor{red}{Additional Constraints on $\beta$} \\

It will be useful to elucidate the relationship between $\beta_{7}$ and dX. Note that we may rewrite the interior product in terms of the Hodge star operator.

\begin{equation}
\beta_{7} = X\lrcorner\varphi = \ast(X\wedge\ast\varphi)
\end{equation}

We may take the Hodge star and then the exterior derivative to get an expression for the co-closure which involves dX. 

\begin{equation}
d\ast\beta_{7} = dX\wedge\ast\varphi - X\wedge\tau\wedge\varphi
\end{equation}

Wedging dX into $\ast\varphi$ isolates $\pi_{7}(dX)$, and we may also simplify the second term into a contraction of X into the torsion two-form.

\begin{equation}
d\ast\beta_{7} = \pi_{7}(dX) \wedge \ast \varphi - \ast(X\lrcorner\tau)
\end{equation}

We know even without assuming compactness that the two-form 2$\beta_{7}$ - $\beta_{14}$ is co-closed. Therefore, if $\beta_{7}$ is also co-closed, then $\beta_{14}$ and $\beta$ itself must be co-closed. This constraint shows how dX and $X\lrcorner\tau$ must be related for $\beta_{7}$ to be co-closed. Assuming compactness, we can see that if $\pi_{7}(dX)$ vanishes then dX must vanish identically. This result follows from an expression for the components of two-forms and an application of Stokes' theorem. 

\begin{equation}
dX\wedge dX\wedge \varphi =  (2 |\pi_{7}(dX)|^{2} - |\pi_{14}(dX)|^{2})\textnormal{vol}
\end{equation}

If dX equals zero, then the constraint shows how the co-closure of $\beta_{7}$ and $X\lrcorner\tau$ must be related.

\section{Exact G$_2$-Structures Subjected to Additional Constraints}

An important subclass of closed G$_{2}$-Structures are the Extremally Ricci-Pinched (ERP) G$_{2}$-Structures. These G$_{2}$-Structures satisfy a condition on the derivative of their torsion. 

\begin{equation}
d\tau = 1/6 (|\tau|^2\varphi + \ast (\tau \wedge \tau))
\end{equation}

Bryant showed that all closed G$_{2}$-Structures which are not torsion-free obey a pinching inequality relating the squared norm of the traceless Ricci tensor to the squared scalar curvature.

\begin{equation}
\int_{M} |\textnormal{Ric}_{0}|^2 \textnormal{vol} \geq 4/21 \int_{M} R^2 \textnormal{vol}
\end{equation}

The ERP G$_{2}$-Structures are those for which this inequality is actually an equality. In this sense ERP G$_{2}$-Structures are special because they are as Ricci-pinched as a G$_{2}$-Structure can be without being torsion-free. \\

ERP G$_2$-Structures actually belong to a larger class of G$_2$-Structures that are called \textit{$\lambda$-quadratic}. The derivative of torsion of these closed G$_2$-Structures depends quadratically on the torsion tensor and features a constant $\lambda$.

\begin{equation}
d\tau = \frac{1}{7} |\tau|^2\varphi + \lambda(\frac{1}{7}\varphi + \ast (\tau \wedge \tau))
\end{equation}

The case $\lambda=1/6$ corresponds to the ERP G$_2$-Structures, but this class of structures also includes familiar closed G$_2$-Structures such as Einstein structures as special cases \cite{Bryant1}. Recently, Ball has shown that many more values for $\lambda$ than were previously understood may be found on non-compact manifolds \cite{Ball}. However, only the ERP condition can occur on a compact manifold so we will focus on that case \cite{Bryant1}. \\

For an exact G$_{2}$-Structure on a compact manifold, we know it cannot be torsion-free because by Hodge theory, if this was the case then it would vanish identically. It still stands that it could be ERP; however, perhaps not suprisingly, we can prove that this cannot happen. 

\begin{theorem} [ERP and Exact G$_{2}$-Structures] \label{thm:ERP} An Exact G$_{2}$-Structure on a compact manifold cannot be Extremally Ricci Pinched. 
\end{theorem}

\textit{Proof}. The ERP condition may be given in terms of the exterior derivatve of the torsion two-form. 

\begin{equation}
d\tau = 1/6 (|\tau|^2\varphi + \ast (\tau \wedge \tau))
\end{equation}

Here $\varphi=d\beta$ as the G$_{2}$-Structure is assumed to be exact. An important property of ERP G$_{2}$-Structures on compact manifolds is that the norm of the torsion two-form is a constant \cite{Bryant1}. As such we may commute this norm with the exterior derivative in the ERP equation.

\begin{equation}
d(\tau - 1/6|\tau|^2\beta) = 1/6 \ast(\tau\wedge\tau)
\end{equation}

We may then write out an expression for the pointwise norm of $\tau \wedge \tau$.

\begin{equation}
\tau \wedge \tau \wedge \ast (\tau \wedge \tau) = |\tau \wedge \tau|^2 \textnormal{vol}
\end{equation}

We may then substitute our expression for $\ast (\tau \wedge \tau)$ into this expression. 

\begin{equation}
6(\tau \wedge \tau \wedge d(\tau - 1/6|\tau|^2\beta) ) = |\tau \wedge \tau|^2\textnormal{vol}
\end{equation}

We may then note that $\tau \wedge \tau$ is automatically closed for an ERP G$_{2}$-Structure \cite{Bryant1}, so we may again use Stokes' theorem to show that $\tau \wedge \tau$ vanishes. We then use a general result for G$_{2}$-Structures that relates the pointwise norm of a differential two form $\eta \in \Lambda^{2}_{14}(T^{*}M)$ to the norm of $\eta \wedge \eta$ \cite{Bryant1}.

\begin{equation}
|\eta\wedge\eta|^2 = |\eta|^4
\end{equation}

As a result we may conclude that $\tau$ vanishes, so the G$_{2}$-Structure is torsion-free. But this would mean that the G$_{2}$ three-form is harmonic, and we already know it is exact. Then by Hodge theory the three-form must vanish, which would contradict positivity. \qed \\

We can apply several of these results on exact G$_2$-Structures to the study of certain closed Laplacian solitons. A \textit{Laplacian soliton} is a triple ($\varphi$, $\bar{X}$, $\lambda$) consisting of a G$_2$-Structure $\varphi$ whose Hodge Laplacian is equal to a constant $\lambda$ times $\varphi$ plus the Lie derivative of $\varphi$ with respect to a vector field $\bar{X}$ .

\begin{equation}
\Delta\varphi = \lambda\varphi + L_{\bar{X}}\varphi
\end{equation}

If ($\varphi_{0},\bar{X},\lambda$) satisfies this equation, then we can use this initial data to write down a solution of the Laplacian flow that essentially scales the initial G$_2$ three-form up to diffeomorphism. The precise expression for $\varphi(t)$ depends on a function $\rho(t)$, which is linear in t, and the action of the group of diffeomorphisms on the initial G$_{2}$-Structure $\varphi_{0}$ which are generated by the vector field $\bar{X}$. 

\begin{equation}
\varphi(t) = \rho(t) \phi_{t}^{*}\varphi_{0}
\end{equation}

If the constant $\lambda$ is positive, negative, or zero, the Laplacian soliton is called \textit{expanding}, \textit{shrinking}, or \textit{steady}, resectively. Lin proves that a Laplacian soliton cannot be shrinking on a compact manifold, and that if it is steady then the G$_2$-Structure is torsion-free \cite{Lin}.  If the G$_2$-Structure is closed, then a Laplacian soliton is called a \textit{closed Laplacian Soliton}. For a closed Laplacian soliton, by results of Lotay and Wei it is known that $\bar{X}$ cannot vanish unless the soliton is steady and the G$_2$-Structure is torsion-free \cite{LW}. They also provide a simple proof that a closed Laplacian soliton cannot be shrinking on a compact manifold and that if it is steady the G$_2$-Structure is torsion-free.  \\

For a compact manifold, then, we will focus on the case $\lambda > 0$. To date, no examples of such expanding closed Laplacian solitons are known to exist on compact manifolds. \\

%Add comment concerning generality of laplacian solitons in space of exact G$_{2}$-Structures

Manipulating the Laplacian soliton equation a bit, it follows that expanding (or shrinking) closed Laplacian solitons are automatically exact G$_2$-Structures. To see this we can express the Hodge Laplacian of the G$_2$-three-form in terms of the derivative of the torsion form.

\begin{equation}
\Delta\varphi = dd^{\ast}\varphi + d^{\ast}d\varphi = dd^{\ast}\varphi = d\tau
\end{equation}

$L_{\bar{X}}\varphi$ also turns out to be an exact form as a consequence of Cartan's formula.

\begin{equation}
L_{\bar{X}}\varphi = d(\bar{X}\lrcorner\varphi) + \bar{X}\lrcorner d\varphi =  d(\bar{X}\lrcorner\varphi)
\end{equation}

We can then rearrange the Laplacian soliton equation to get an expression for $\lambda\varphi$. Note that by assumption the constant $\lambda$ doesn't vanish. 

\begin{equation}
\lambda\varphi = d(\tau - \bar{X}\lrcorner\varphi)
\end{equation}

As a consequence, all expanding closed Laplacian solitons are actually exact, so if exact G$_2$-Structures generally can't exist on compact manifolds, then a closed Laplacian soliton would have to have $\lambda=0$ and the G$_2$-Structure would have to be torsion-free. \\

Not only does the Laplacian soliton equation imply that $\varphi$ is exact, but it also gives information about the two-form. Fixing a choice of gauge, we can take the two-form to be exactly $\lambda\beta=\tau - \bar{X}\lrcorner\varphi$. Given this choice, the decomposition of $\beta$ into irreducible G$_2$ representations is immediate because each of these terms are in particular G$_2$ representations. 

\begin{align}
&\lambda\beta_{7} = \lambda X\lrcorner\varphi = -\bar{X}\lrcorner\varphi \\
&\lambda\beta_{14} = \tau
\end{align}

We can now apply some of the results we found for general exact G$_2$-Structures to expanding closed Laplacian solitons. It turns out that several of the key results found by Lotay and Wei are immediate consequences of the fact that the G$_{2}$-Structure underlying (non-steady) closed Laplacian solitons is exact. 

\begin{theorem} [Projection Identities for Shrinking and Expanding Closed Laplacian Solitons] \label{thm: LaplacianSolitonProjection} Given a Closed Laplacian Soliton $(\varphi, \bar{X},\lambda)$ with $\lambda \neq 0$, the norm of the torsion, the constant $\lambda$, and $d^{\ast}\bar{X}$ satisfy a linear relation. Moreover, the form $\bar{X}\lrcorner\varphi$ is co-closed.
\end{theorem}

\begin{align}
& \frac{1}{3} |\tau|^2 + d^{\ast}\bar{X} = 7/3 \lambda \\
& d^{\ast} (\bar{X}\lrcorner\varphi) = 0
\end{align}

\textit{Proof.} The first relation is just the statement $\pi_{1}(d\beta_{7}) + \pi_{1}(d\beta_{14}) = 7$ applied to the case of a Laplacian soliton. For such a structure, the components of the exact G$_2$ two-form are given in terms of $\tau$ and $\bar{X}\lrcorner\varphi$. We saw that we may simplify this expression into an expression relating the co-exterior derivative of the vector field X and the pointwise inner product of $\beta_{14}$ and $\tau$.

\begin{equation}
\langle\beta_{14},\tau\rangle + d^{\ast}X = 1
\end{equation}

For the Laplacian soliton, $\beta_{7}$ is given in terms of the vector field $\bar{X}$ and $\beta_{14}$ is given in terms of the torsion two-form. 

\begin{align}
&\lambda\beta_{7} = \lambda X\lrcorner\varphi = -\bar{X}\lrcorner\varphi \\
&\lambda\beta_{14} = \tau
\end{align}

We can then plug these expressions for the components into the general projection expression for an exact G$_2$-Structure to provide the desired relation. \\

We can now look at the relation $\pi_{7}(d\beta_{7}) + \pi_{7}(d\beta_{14}) = 0$ for a Laplacian soliton. We saw as a consequence of this relation that there was a canonical co-closed two-form built out of the components of the exact G$_2$-Structure. 

\begin{equation}
d^{\ast} (2\beta_{7}-\beta_{4}) = 0
\end{equation}

For the Laplacian soliton, not only is this form co-closed, but $\beta_{14}$ is also co-closed since it is directly proportional to $\tau$. As a consequence $\beta_{7}$ also has to be co-closed, which proves the desired relation. \qed \\

% Comment on triviality of result for steady solitons

As a consequence, the two-form (with this choice of gauge) underlying the G$_2$-Structure of a Laplacian soliton is co-closed.\\

%\textcolor{red}{Discuss integrating this} \\

We can also apply our result relating the L$^{2}$ norms of $\beta_{7}$ and $\beta_{14}$ to the case of these closed Laplacian solitons.

\begin{theorem} For an expanding closed Laplacian soliton ($\varphi, \bar{X}, \lambda$) on a compact manifold the L$^{2}$ norms of the torsion $\tau$ and the vector field $\bar{X}$ are proportional up to a factor of six.

\begin{equation}
|\tau|^2 = 6|\bar{X}|^2
\end{equation}
\end{theorem}

\textit{Proof.} We plug in the components $\lambda\beta_{7}=-\bar{X}\lrcorner\varphi$ and $\lambda\beta_{14}=\tau$ into the result from Theorem 4.1. 

\begin{equation}
2|\bar{X}\lrcorner\varphi|_{L^{2}}^2 = |\tau|_{L^{2}}^2
\end{equation}

We can then relate the norm of $\bar{X}\lrcorner\varphi$ to that of $\bar{X}$. It is generally true for vector fields that the norm of $\bar{X}$ is three times the norm of $\bar{X}\lrcorner\varphi$ \cite{Karigiannis1}. The result follows. \qed \\

This result provides an alternate proof of the fact that $\bar{X}$ cannot vanish for a closed Laplacian soliton on a compact manifold unless the G$_{2}$-Structure is torsion-free and the soliton is steady. This precise relationship between the norms of $\tau$ and $\bar{X}$, however, appears to be new. Note that $\lambda$ does not enter into the norm expression; however, we have assumed that it does not vanish because it is only in this case that the G$_{2}$-Structure is exact and the method of proof applies. \\

We have also seen that an exact G$_{2}$-Structure cannot by ERP on a compact manifold. As an ERP G$_{2}$-Structure is not torsion-free and thus cannot correspond to a steady soliton, if an ERP G$_{2}$-Structure did correspond to a soliton the G$_{2}$-Structure would have to be exact. We may then infer that an ERP G$_{2}$-Structure may never be a closed Laplacian soliton on a compact manifold. \\

We may also clarify a constraint on $\bar{X}\lrcorner\tau$ found by Lotay and Wei for gradient Laplacian solitons. For these solitons the dual one-form of the vector field $\bar{X}$ is the gradient of a function $f$. We saw the constraint relating the co-closure of $\beta_{7}$ to $dX$ and $X\lrcorner\tau$, which holds for a general exact G$_{2}$-Structure. We may multiply through by $\lambda$ to rephrase this constraint for a closed Laplacian soliton.

\begin{equation}
d\ast (\bar{X}\lrcorner\varphi) = \pi_{7}(d\bar{X}) \wedge \ast \varphi - \ast(\bar{X}\lrcorner\tau)
\end{equation}

For a closed Laplacian soliton, $\beta_{7}$ is co-closed, so the left-hand side vanishes. For a gradient soliton $d\bar{X}$ also has to vanish, so these together imply that $\bar{X}\lrcorner\tau$ has to vanish, which is the result found by Lotay and Wei. From this result we may also establish a converse statement which provides a means to ensure that a soliton is a gradient soliton.

\begin{theorem} [Gradient Solitons] \label{thm:gradient} A necessary and sufficient criterion for a Laplacian soliton $(\varphi, \bar{X}, \lambda)$ to be a gradient soliton on a simply-connected compact manifold is for $\bar{X}\lrcorner\tau=0$.
\end{theorem}

\textit{Proof.} Lotay and Wei prove that gradient solitons, without an assumption on compactness, satisfy $\bar{X}\lrcorner\tau=0$. For the converse we may utilize the same relation. 

\begin{equation}
d\ast (\bar{X}\lrcorner\varphi) = \pi_{7}(d\bar{X}) \wedge \ast \varphi - \ast(\bar{X}\lrcorner\tau)
\end{equation}

If $\bar{X}\lrcorner\tau$ vanishes then this implies that $\pi_{7}(d\bar{X})$ has to vanish pointwise because $\beta_{7}$ always is co-closed for a Laplacian soliton. However, on a compact manifold the $L^2$ norms of $\pi_{7}(d\bar{X})$ and $\pi_{14}(d\bar{X})$ are related up to a factor of two. 

\begin{equation} 
2 |\pi_{7}(d\bar{X})|_{L^2}^{2} = |\pi_{14}(d\bar{X})|_{L^2}^{2}
\end{equation}

Therefore, $\pi_{14}(d\bar{X})$ also has to vanish so $d\bar{X}$ has to vanish. If the underlying manifold is simply-connected or more generally has vanishing first Betti number, then $\bar{X}=df$ for some function $f$ and the soliton is a gradient soliton. \qed \\

The significance of gradient solitons is derived by analogy to the Ricci flow. In that context Perelman showed that every Ricci soliton is a gradient soliton using variational methods \cite{Perelman}. If the Laplacian flow behaves analogously to the Ricci flow then we may also expect gradient solitons to be generic. However, the precise degree of analogy between these geometric flows is rather unclear. In particular, it is unclear if gradient solitons exist at all on compact manifolds. However, there are examples of complete steady gradient solitons recently constructed by Ball \cite{Ball}. 

\section{The Curvature Tensor of Exact G$_2$-Structures}

In addition to several results on the Ricci curvature for closed G$_{2}$-Structures, there are also several results on the Riemann curvature which were proven by Cleyton and Ivanov \cite{CI1, CI2}. In this section we show how to strengthen these results on the Riemann curvature when the G$_{2}$-Structure is exact. \\

Recall that for a general n-dimensional Riemannian manifold, we may decompose its Riemann curvature as a (0,4) tensor orthogonally into components \cite{Besse}.

\begin{equation}
R_{ijkl} = S_{ijkl} + E_{ijkl} + W_{ijkl}
\end{equation}

These components may be expressed in terms of the metric, scalar curvature, and traceless Ricci tensor Ric$_{0}$. 

\begin{align}
& S_{ijkl} = \frac{R}{n(n-1)} (g_{il}g_{jk}-g_{ik}g_{jl}) \\
& E_{ijkl} = \frac{1}{n-2} ((\textnormal{Ric}_{0})_{il}g_{jk} - (\textnormal{Ric}_{0})_{jl}g_{ik} - (\textnormal{Ric}_{0})_{ik}g_{jl} + (\textnormal{Ric}_{0})_{jk}g_{il}) \\
& W_{ijkl} = R_{ijkl} -S_{ijkl} - E_{ijkl}
\end{align} 

These tensors have the same algebraic symmetries as the Riemann tensor. The Weyl tensor W is also traceless and is known to measure the deviation of the manifold from being locally conformally flat \cite{Besse}.  \\

The previous decomposition may be thought of as a decomposition of the space of algebraic curvature tensors into O(n) representations. In the G$_{2}$ context we can decompose the space of algebraic curvature tensor further into irreducible G$_2$ representations. This decomposition was carried out by Cleyton and Ivanov, who found that the Weyl tensor from the SO(7) decomposition further splits into three pieces. 

\begin{equation}
W = W_{27} + W_{64} + W_{77}
\end{equation}

This full decomposition (with the terms S and E as well) shows how to decompose the space of algebraic curvature tensors into irreducible representations of G$_{2}$. The Weyl tensor spaces found by Cleyton and Ivanov are isomorphic to familiar G$_{2}$ representations.

\begin{align}
& W_{77} \cong V_{77}^{(2,0)} \\
& W_{64} \cong V_{64}^{(1,1)} \\
& W_{27} \cong \Lambda^{3}_{27}
\end{align}

These spaces may be defined explicitly in terms of the Ricci curvature and scalar curvature of the G$_{2}$-Structure, as well in terms of as a tensor that Cleyton and Ivanov refer to as the $\phi$-Ricci Tensor \cite{CI1}. We will not need the precise definitions of these spaces for our purposes but it is useful to keep in mind the representations of G$_{2}$ to which they are isomorphic. \\

In order to prove results on the curvature, we study the first Pontryagin class of the underlying manifold. Given an affine connection $\nabla$ we can build a representative of this class using the curvature of the connection. 

\begin{equation}
p_{1}(\nabla) = -\frac{1}{8\pi^2} \textnormal{tr}(R_{\nabla}^2)
\end{equation}

In their proof that a closed G$_{2}$-Structure on a compact manifold cannot be Einstein unless it is torsion-free, Cleyton and Ivanov consider the integral of the closed G$_{2}$-Structure wedged against the difference of the Pontryagin forms corresponding to the Levi-Civita connection and the canonical G$_{2}$ connection. 

\begin{equation}
\int_{M}( p_{1}(\nabla_{LC})-p_{1}(\nabla_{G_{2}})) \wedge \varphi
\end{equation}

As the class is topological and in particular doesn't depend on the choice of connection, the difference of the representatives of the Pontryagin classes of the two connections is an exact four-form. Since $\varphi$ is closed, we can use Stokes' theorem to assert that this integral vanishes. If we assume the stronger condition that the G$_{2}$-Structure is exact then we do not need to take the difference of two Pontryagin class representatives to utilize Stokes' theorem. Since the Pontryagin form for a given connection is a closed form, we may immediately say that the integral of such a form wedged against the G$_{2}$ three-form vanishes. 

\begin{equation}
\int_{M} p_{1}(\nabla) \wedge \varphi = 0 \\
\end{equation}

We may then rewrite these identities in terms of the curvatures of the connections to get constraints on the curvature of an exact G$_{2}$-Structure. The Levi-Civita constraint is particularly easy to work with due to the symmetries of the Riemann curvature, and immediately provides a result on the G$_2$ decomposition of the Riemann cuvature viewed as a Lie algebra-valued two-form. 

\begin{theorem} [Integral Identity for two-form components of Riemann curvature of Exact G$_{2}$-Structure] \label{thm:CurvatureComponents1} If Rm is the Riemann curvature of the metric determined by an exact G$_{2}$-Structure on a compact manifold then the L$^{2}$ norms of the irreducible G$_2$ components of the Riemann curvature are related up to a factor of two. 
\end{theorem}

\begin{equation}
2 | \pi_{7}(\textnormal{Rm}) |_{L^{2}}^2 = | \pi_{14}(\textnormal{Rm}) |_{L^{2}}^2 
\end{equation}

\textit{Proof}. We start with the Pontryagin class identity for the Levi-Civita connection, substituting in the expression for the representative of the class in terms of the curvature. 

\begin{equation}
\int_{M} \textnormal{tr}(\textnormal{Rm} \wedge \textnormal{Rm}) \wedge \varphi =0
\end{equation}

This expression is reminiscent of an expression for the components of an ordinary two-form $\alpha$ that we have utilized in several preceding theorems.

\begin{equation}
\alpha \wedge \alpha \wedge \varphi = 2|\alpha_{7}|^2 - |\alpha_{14}|^2
\end{equation}

Viewing the curvature as a Lie algebra-valued two-form and using the symmetries of the Riemann curvature, we get an analogous identity. The result then follows. \qed \\

We may compare this result to what Joyce finds for a torsion-free G$_{2}$-Structure. For such a structure, the curvature only has a component in the 14 dimensional representation of G$_{2}$ on two-forms. As a result, for a torsion-free G$_2$-Structure, wedging the G$_2$ three-form against the Levi-Civita representative of the first Pontryagin class and integrating gives the L$^{2}$ norm of the curvature. 

\begin{equation}
\int_{M} p_{1}(\nabla_{LC}) \wedge \varphi = -\int_{M} |R|^2 \textnormal{vol}
\end{equation}

From this expression we see that the first Pontryagin class being trivial would imply that the connection is flat; in other words, that it would have trivial holonomy. However, for an exact G$_{2}$-Structure the integral always vanishes. This is consistent with a vanishing first Pontryagin class but not necessary. So, at least by this reasoning, the first Pontryagin class vanishing is not an obstruction for a compact manifold to admit an exact G$_{2}$-Structure. \\

We may also point out that it is known for an arbitrary G$_2$-Structure that $\pi_{7}(\textnormal{Rm})$ cannot vanish unless the G$_2$-Structure is torsion-free \cite{Karigiannis2}. Moreover, if $\pi_{14}(\textnormal{Rm})$ vanishes then the G$_{2}$-Structure is not only torsion-free but flat. So we could have anticipated that neither of these curvature components could vanish for an exact G$_2$-Structure, but it is certainly not true for an arbitrary G$_2$-Structure that these components would be related precisely as in this theorem. \\

We may also express the wedge product of the Pontryagin form against the G$_{2}$-Structure in terms of the decomposition of the Riemann tensor due to Cleyton and Ivanov. They find a relationship between the G$_2$ components of the Riemann tensor and the first Pontryagin class for an arbitrary closed G$_2$-Structure on a compact manifold \cite{CI1}. 

\begin{equation}
\langle p_{1}(M) \cup \varphi, [M]\rangle \geq -\frac{1}{8\pi^2} \int_{M} |W_{77}|^2 - 1/2|W_{64}|^2 -3/16 S^2
\end{equation} 

Here equality holds iff the closed G$_2$-Structure is ERP. We can now consider the implications of this identity when the G$_2$-Structure is exact. 

\begin{theorem} [Integral Identity for components of Riemann curvature of Exact G$_{2}$-Structure] \label{thm:CurvatureComponents2} For an exact G$_2$-Structure on a compact manifold M, the L$^{2}$ norms of the components of the Riemann tensor are related by an inequality. 
\end{theorem}

\begin{equation}
 \int_{M}|W_{77}|^2 > 1/2 \int_{M}|W_{64}|^2 + 3/16 \int_{M}S^2
\end{equation}

\textit{Proof}. We start with the identity of Cleyton and Ivanov and note that for an exact G$_2$-Structure, the integral of the first Pontryagin class wedged against the G$_2$ three-form vanishes.

\begin{equation}
0 \geq  -\int_{M}|W_{77}|^2 + 1/2\int_{M}|W_{64}|^2 + 3/16 \int_{M} S^2
\end{equation}

Rearranging, we basically get the required result, but we can also make the inequality a proper inequality by noting that the case of equality only happens when the closed G$_2$-Structure is ERP. We proved in theorem 4.2 that an exact G$_2$-Structure cannot be ERP, so this statement cannot be an equality. \qed \\

As a corollary we may see that W$_{77}$ may not vanish for an exact G$_2$-Structure. This result complements the finding by Cleyton and Ivanov that W$_{27}$ also cannot vanish for a closed G$_{2}$-Structure unless that structure is torsion-free. Of course, S$^2$ also may not vanish since that would also mean that the G$_2$-Structure was torsion-free. \\

%\textcolor{red}{Comment on 64 being related to ERP} \\

\subsection{Acknowledgements} I would like to thank Gavin Ball, Robert Bryant, Jason Lotay, and my advisor Dave Morrison for helpful comments related to my work on exact G$_{2}$-Structures. I'm also very appreciative of the support I have recieved more generally from the Simons Collaboration on Special Holonomy in Geometry, Analysis, and Physics. I'd also like to thank my friends and family for support and Jessica Li in particular for being very critical of my grammar which resulted in noticeable improvements for this paper.

 \end{document}